\input amstex\documentstyle{amsppt}  
\pagewidth{12.5cm}\pageheight{19cm}\magnification\magstep1
\topmatter
\title A monotonicity property for the new basis of
$\bold C[(\bold Z/2)^D]$
\endtitle
\author G. Lusztig\endauthor
\address{Department of Mathematics, M.I.T., Cambridge, MA 02139}\endaddress
\thanks{Supported by NSF grant DMS-2153741}\endthanks
\endtopmatter   
\document
\define\mat{\matrix}
\define\endmat{\endmatrix}

\define\bx{\boxed}

\define\pe{\perp}
\define\si{\sim}

\define\sqc{\sqcup}

\define\op{\oplus}
   
\define\part{\partial}
\define\emp{\emptyset}

\define\n{\notin}

\define\m{\mapsto}
\define\do{\dots}

\define\sub{\subset}    

\define\T{\times}

\define\nl{\newline}
\redefine\i{^{-1}}

\redefine\spa{\spadesuit}

\define\e{\epsilon}

\define\ph{\phi}
\define\ps{\psi}

\redefine\t{\tau}
\define\th{\theta}

\define\CC{\bold C}

\define\NN{\bold N}

\define\ZZ{\bold Z}

\define\cf{\Cal F}

\define\ci{\Cal I}

\head Introduction \endhead
\subhead 0.1\endsubhead
Let $V_D$ be a vector space of finite even dimension $D\ge2$
over a field $F$ endowed with a nondegenerate symplectic form
$(,):V_D\T V_D@>>>F$ and with a collection of vectors
$e_1,e_2,\do,e_D,e_{D+1}$ such that 

$(e_i,e_j)=1$ if $j=i+1\mod D+1$, $(e_i,e_j)=-1$ if
$i=j+1\mod D+1$ and $<e_i,e_j>=0$ for all other pairs $i,j$.
\nl
Such a collection is called a ``circular basis''.
(It always exists and satisfies $\sum_{i=1}^{D+1} e_i=0$.)

Let $\CC[V_D]$ be the $\CC$-vector space with basis
$\{\bx{x};x\in V_D\}$. 

\subhead 0.2\endsubhead
In the rest of the paper (except in 3.1 and 3.7) we assume that
$F$ is the field with two elements. 
In this case the definition of circular basis appeared in \cite{L20a}.

In \cite{L20a},\cite{L23} a new basis of $\CC[V_D]$ was defined
in terms of a circular basis of $V_D$. (This new basis was in fact
defined earlier in \cite{L20} in a less symmetric way.)
The elements of the new basis are characteristic functions of certain isotropic
subspaces of $V_D$. They can be expressed
in terms of the old basis by a matrix which is upper triangular with
respect to a certain partial order $\le$ on the new basis. In this
paper we prove a monotonicity property of the new basis, namely that if
$B,B'$ in the new basis satisfy $B\le B'$ then the dimensions $d,d'$
of the corresponding isotropic subspaces satisfy $d\le d'$.
The proof is based on a construction given in \S1 which
extends and simplifies results in \cite{L20, 1.18}.

In \S3 we state a conjecture (unrelated to the monotonicity property)
on the entries of the Fourier transform $\CC[V_D]@>>>\CC[V_D]$
with respect to the new basis.

\subhead 0.3\endsubhead
{\it Notation} For $i,j$ in $\ZZ$ we set $[i,j]=\{z\in\ZZ;i\le z\le j\}$.
For a finite set $?$ we denote by $|?|$ the cardinal of $?$.

\head 1. The assignment $B\m B[i]$\endhead
\subhead 1.1\endsubhead
For any $I\sub[1,D+1]$ we set $e_I=\sum_{s\in I}e_s\in V_D$.

Let $\ci$ be the set of all $I\sub[1,D+1]$ of the
form $[i,j]$ or $[1,D+1]-[i,j]$ with $1\le i\le j\le D$. (Such sets are
called intervals.)
We have $\ci=\ci^0\sqc\ci^1$
where $\ci^0=\{I\in\ci;|I|=0\mod2\}$, $\ci^1=\{I\in\ci;|I|=1\mod2\}$.
For $I,I'$ in $\ci^1$ we write $I\prec I'$ whenever
$I\subsetneqq I'$ and $I'-I=I_1\sqc I_2$ with $I_1\in\ci,I_2\in\ci$.
For $I,I'$ in $\ci^1$ we write $I\spa I'$ whenever $I\cap I'=\emp$
and $I\cup I'\n\ci$.
For $I\in\ci^1$ let $I^{ev}$ be the set of all $i\in I$ such that
$I-\{i\}=I'\sqc I''$, with $I'\in\ci^1$, $I''\in\ci^1$, $I'\spa I''$.
Let $I^{odd}=I-I^{ev}$. We have $|I^{ev}|=(|I|-1)/2$. 

Let $R$ be the set whose elements are finite unordered sequences of
objects of $\ci^1$.

For $B\in R$ let $<B>$ be the subspace of $V_D$
generated by $\{e_I;I\in B\}$.

For $i\in[1,D+1]$, $B\in R$ we set $B_i=\{I\in B;i\in I\}$, 
$$g_i(B)=|B_i|$$
and
$$\e_i(B)=(1/2)g_i(B)(g_i(B)+1)\in F.$$
For $B\in R$ we set
$$\e(B)=\sum_{i\in[1,D+1]}\e_i(B)e_i\in V_D.$$

Following \cite{L23} we define
$\ph(V_D)$ to be the set consisting of all $B\in R$ such that
$(P_0),(P_1)$ below hold.

$(P_0)$ If $I\in B,I'\in B$, then $I=I'$, or $I\spa I'$, or
$I\prec I'$, or $I'\prec I$.

$(P_1)$ Let $I\in B$. There exist $I_1,I_2,\do,I_k$ in $B$ such that
$I^{ev}\sub I_1\cup I_2\cup\do\cup I_k$ (disjoint union),
$I_1\prec I$, $I_2\prec I,\do,I_k\prec I$.
\nl
For $B',B$ in $\ph(V_D)$ we say that $B'\le B$ if there exist
$B^0,B^1,B^2,\do,B^k$ in $\ph(V_D)$ such that $B^0=B',B^k=B$,

$\e(B^0)\in<B^1>,\e(B^1)\in<B^2>,\do,\e(B^{k-1})\in<B^k>$.
\nl
By \cite{L20},\cite{L23},

(a) $\le$ is a partial order on $\ph(V_D)$;

(b) $\e:\ph(V_D)@>>>V_D$ is a bijection and $\e(B)\in<B>$ for any
$B\in\ph(V_D)$;

(c) For $B\in\ph(V_D)$, $\{e_I;I\in B\}$ is a basis of $<B>$; hence
$\dim<B>=|B|$.
\nl
For $x\in V_D$ we define $B(x)\in\ph(V_D)$ by $\e(B(x))=x$.
For $x,y$ in $V_D$ we write $x\le y$ whenever $B(x)\le B(y)$. By (a),
this is a partial order on $V_D$; it follows that

(d) $\{\bx{<B>};B\in\ph(V_D)\}$ is a ``new'' basis of $\CC[V_D]$
related to the
basis $\{\bx{x};x\in V_D\}$ by an upper triangular matrix with
$1$ on diagonal.

Let $x\in V_D$. We have $0\le x$ for any
$x\in V_D$. Hence the set of sequences

$\{0=x_0,x_1,\do,x_k=x\text{ in $V_D$ such that }x_0<x_1<\do<x_k\}$
\nl
is nonempty and we must
have $k\le|V_D|$ (since $\le$ is a partial order). Hence we can
define $\nu(x)\in\NN$ as the maximum value of $k$ over
all sequences as above.

\subhead 1.2\endsubhead
Let $B\in\ph(V_D),i\in[1,D+1]$ be such that $\{i\}\in B$. We can
arrange the elements of $B_i$ in a sequence $I_1,I_2,\do,I_k$ with
$\{i\}=I_1\prec I_2\prec\do\prec I_k$, $k=g_i(B)\ge1$.
We show:

(a) {\it If $k\ge2$, then $i\in I_2^{ev}$.}
\nl
Assume that $i\in I_2^{odd}$. Then we have either $I_2-\{i\}\in\ci$
(contradicting $\{i\}\prec I_2$) or $I_2-\{i\}=I'\sqc I''$ where
$I'\in\ci^0, I''\in\ci^0$. In this last case we set $i'=i+1$ (if
$i\in[1,D]$), $i'=1$ (if $i=D+1$). We have $i'\in I'$ or $i'\in I''$.
We can assume that $i'\in I'$. Since $i\in I_2^{odd}$ we have
$i'\in I_2^{ev}$. Hence there exists $J\in B$ such that $J\prec I_2$.
$i'\in J$. By the choice of $i'$ we have $J\cup\{i\}\in\ci$.
If $i\n J$ then by (P1) we have $\{i\}\spa J$ so that
$J\cup\{i\}\n\ci$; this is a contradiction. We see that $i\in J$.
Since $J\in B_i$ and $J\prec I_2$, we must have $J=\{i\}$. This
contradicts $i'\in J$; (a) is proved.

We show:

(b) {\it If $k\ge3$, then $i\in I_3^{odd}$.}
\nl
We can find  $J_1,J_2,\do,J_s$ in $B$ such that
$I_3^{ev}\sub J_1\cup J_2\cup\do\cup J_s$ (disjoint union),
$J_1\prec I_3$, $J_2\prec I_3,\do,J_s\prec I_3$.
From the definitions we have

(c) $J_1^{ev}\cup J_2^{ev}\cup\do\cup J_s^{ev}\sub I_3^{odd}$.
\nl
We have $i\in J_1\cup J_2\cup\do\cup J_s$. We can assume that
$i\in J_1$ so that $J_1\in B_i$. Since $J_1\prec I_3$ we have
$J_1=I_2$ or $J_1=I_1=\{i\}$.

Assume first that $J_1=I_2$. Since $i\in I_2^{ev}=J_1^{ev}$
(see (a)), we see from (c) that $i\in I_3^{odd}$, as desired.

In the rest of the proof we assume that $J_1=I_1=\{i\}$.
Since $\{i\}\prec I_2$ we have $|I_2|\ge3$.
Hence there are well defined elements $j\ne j'$ in $I_2$
such that $I_2-\{j\}\in\ci^0,I_2-\{j'\}\in\ci^0$.

Assume now that $j\in I_3^{odd}$. 
Since $I_2\prec I_3$ we can find $h\in I_3^{ev}$
such that $h\n I_2$, $\{h,j\}\in\ci^0$.
By $(P_1)$ we can find $t\in[1,s]$ such that $J_t\prec I_3$, $h\in J_t$.
Since $h\in J_t,h\n I_2$ we see that $J_t\not\sub I_2$.
We have $I_2\not\sub J_t$. (If $I_2\sub J_t$, we would have
$J_1=I_1\sub J_t$
hence $t=1$ so that $J_t=\{i\}$ and $h=i$ contradicting $h\n I_2$).
From $(P_0)$ we now see that $J_t\spa I_2$. But this contradicts
$h\in J_t,\{h,j\}\in\ci^0,j\in I_2$.

We see that $j\in I_3^{ev}$. Similarly we have $j'\in I_3^{ev}$.
Since $i\in I_2^{ev}$, this implies that $i\in I_3^{odd}$.
This completes the proof of (b).

\subhead 1.3\endsubhead
Let $B\in\ph(V),i\in[1,D+1]$ be such that $\{i\}\in B$.
We arrange the elements of $B_i$ in a sequence $I_1,I_2,\do,I_k$
as in 1.2. Recall that $k=g_i(B)\ge1$. We define $B[i]\in R$ as follows.
If $g_i(B)=1$ then $B[i]=B-\{i\}$.
If $g_i(B)\ge2$ then $B[i]=(B\sqc\{H\}\sqc\{H'\})-(\{I_1\}\sqc\{I_2\})$
where $H,H'$ in $\ci^1$ are defined by $I_2-\{i\}=H\sqc H'$.
\nl
(Here we use that $i\in I_2^{ev}$, see 1.2(a) and that $H\n B,H'\n B$
since neither $H,I_1$ or $H',I_1$ satisfy $(P_0)$.) We show:

(a) $B[i]\in\ph(V)$.
\nl
We first show:

(b) {\it any $I,I'$ in $B[i]$ satisfy $(P_0)$.}
\nl
If $g_i(B)=1$ then (b) is obvious. We now assume that $g_i(B)\ge2$.
If $I$ and $I'$ are not in $\{H,H'\}$ then (b) is obvious.
Assume that $I\n\{H,H'\},I'\in\{H,H'\}$. Then $I\ne I'$.
If $I\spa I_2$ then $I\spa I'$. If $I_2\prec I$ then $I'\prec I$.
If $I\prec I_2$ and $i\n I$ then $I\sub H$ or $I\sub H'$; moreover
if $I\sub H$ we have $I\prec H$ unless $I\cup\{i\}\in\ci^0$ (but this
would contradict $(P_0)$ for $I,I_1$). Thus if
$I\prec I_2$ and $i\n I$ then $I\prec H$ or $\prec H'$.
If $I\prec I_2$ and $i\in I$ then $I=I_1$ so that $I\n B[i]$.
If $I_2\spa I$ then $I\spa H$ and $I\spa H'$. We see that in our case
(b) is satisfied.
The same argument applies when $I\in\{H,H'\},I'\n\{H,H'\}$. 
Assume now that $I\in\{H,H'\},I'\in\{H,H'\}$. Then we have either
$I=I'$ or else $I\spa I'$. This proves (b).

We show:

(c) {\it any $I\in B[i]$ satisfies $(P_1)$.}
\nl
Assume first that $g_i(B)=1$. We have $I\in B,I\ne\{i\}$.
Let $J_1,J_2,\do,J_k$
be as in $(P_1)$ for $I$ (relative to $B$). Since $g_i(B)=1$
and we have $i\n I$ hence $i\n J_1,i\n J_2,\do,i\n J_k$. Thus
$J_1,J_2,\do,J_k$ are as in $(P_1)$ for $I$ (relative to $B[i]$).
This proves (c) in our case.

Next we assume that $g_i(B)\ge2$. Assume that $I\in B$,
$I\ne I_2,I\ne I_1$. Let $J_1,J_2,\do,J_k$ be as in $(P_1)$ for $I$
(relative to $B$). If all $J_s$ are different from $I_2,I_1$ then
$J_1,J_2,\do,J_k$ are as in $(P_1)$ for $I$ (relative to $B[i]$).
If one $J_s$ (say $J_1$) is equal to $I_2$ (so that $J_2,\do,J_k$
are different from $I_1$ since $I_1\sub I_2$) then we must have
$I_2\prec I$ hence $I$ is one of $I_3,I_4,\do$. If $i\in I^{odd}$
then the collection $H,H',J_2,J_3,\do,J_k$ is as in $(P_1)$ for $I$
(relative to $B[i]$). If $i\in I^{ev}$ (so that, by 1.2(b),
we have $I\ne I_3$ that is $I$ is one of $I_4,I_5,\do$), then the
collection $I_3,J_2,J_3,\do$ satisfies the conditions of $(P_1)$
for $I$ (relative to $B[i]$) except that the condition of ``disjoint
support'' may not be satisfied. But if from this collection we remove
those $J_s$ which are contained in $I_3$, the resulting collection
satisfies the conditions of $(P_1)$ for $I$ (relative to $B[i]$).

Next we assume that $I=H$.
Let $J_1,J_2,\do,J_k$ be as in $(P_1)$ for $I_2$ (relative to $B$).
Let $j\in H^{ev}$. Then $j\in I_2^{ev}$ hence we can find $s$ such that
$j\in J_s^{odd}$. If $i\in J_s$, then $J_s$ must be $\{i\}$ so that
$j=i$, a contradiction. Thus we have $i\n J_s$. It follows that
$J_s\sub H$ or $J_s\sub H'$. Since $j\in H\cap J_s$, we see that
$J_s\sub H$. Using this and the inclusions $j\in J_s^{odd},j\in H^{ev}$,
we deduce $J_s\prec H$.
We have $J_s\ne\{i\}$ since $i\n H$ and $J_s\ne I_2$ since
$J_s\prec I_2$. Thus $J_s\in B[i]$. We see that
$H^{ev}$ is contained in the union of a subcollection
of $J_1,J_2,\do,J_k$ which shows that $(P_1)$ holds for $H$ (relative
to $B[i]$. The same argument applies to $I=H'$.
This completes the proof of (c).

\subhead 1.4\endsubhead
Let $B\in\ph(V_D),i\in[1,D+1]$ be such that $\{i\}\in B$.

If $g_i(B)=1$ then $g_i(B[i])=0$ hence
$$(1/2)g_i(B[i])(g_i(B[i])+1)=(1/2)g_i(B)(g_i(B)+1)+1\mod2.$$
(that is $0=1+1\mod2$).

If $g_i(B)\ge2$ then $g_i(B[i])=g_i(B)-2$ hence
$$(1/2)g_i(B[i])(g_i(B[i])+1)=(1/2)g_i(B)(g_i(B)+1)+1\mod2.$$

If $j\in[1,D+1]-\{i\}$ then $g_j(B[i])=g_j(B)$ hence
$$(1/2)g_j(B[i])(g_j(B[i])+1)=(1/2)g_j(B)(g_j(B)+1).$$

It follows that

$\e(B[i])=\e(B)+e_i\in F$.

\subhead 1.5\endsubhead
Let $B\in\ph(V_D),i\in[1,D+1]$.
The following two conditions are equivalent:
(1) $\{i\}\in B$.
(2) $g_i(B)\ge1$ and $g_{i-1}(B)=g_{i+1}(B)=g_i(B)-1$.
(When $i=1$ we interpret $g_0(B)$ as $g_{D+1}(B)$; when $i=D+1$
we interpret $g_{D+2}(B)$ as $g_1(B)$.)

\head 2. A monotonicity property\endhead
\subhead 2.1\endsubhead
We now assume that $D\ge4$. Let
$$V_{D-2},(,):V_{D-2}\T V_{D-2}@>>>F,e'_1,\do,e'_{D-2},e'_{D-1}$$
be the analogues of $V_D,(,),e_1,\do,e_D,e_{D+1}$
in 1.1 when $D$ is replaced by $D-2$.
For any $j\in[1,D+1]$ let $\t_j:V_{D-2}@>>>V_D$ be the linear map
which takes $e'_1,\do,e'_{D-2},e'_{D-1}$ (in the order written) to:

$e_3,e_4,\do,e_D,e_{D+1}+e_1+e_2$ if $j=1$,

$e_1,\do,e_{j-2},e_{j-1}+e_j+e_{j+1},e_{j+2},\do,e_{D+1}$ if
$j=2,3,\do,D$,

$e_D+e_{D+1}+e_1,e_2,e_3,\do,e_{D-1}$ if $j=D+1$.

This map is injective and compatible with $(,)$.
Its image is a complement of the line $Fe_j$ in $\{x\in V_D;(x,e_j)=0\}$.

For any $E\sub V_D$ we set $\bx{E}=\sum_{x\in E}\bx{x}\in\CC[V_D]$.
For $E'\sub V_{D-2}$ we define $\bx{E'}\in\CC[V_{D-2}]$ in a similar way.

Let $\th_j:\CC[V_{D-2}]@>>>\CC[V_D]$ be the linear map defined by
$\th_j(\bx{x'})=\bx{\t_j(x')}+\bx{\t_j(x')+e_j}$ for any $x'\in V_{D-2}$.

\subhead 2.2\endsubhead
Following \cite{L20a} we define
a collection $\ph'(V_D)$ of subspaces of $V_D$ by induction on $D$.

If $D=2$, $\ph'(V_D)$ consists of the $4$ subspaces of
dimension $\le1$ of $V_D$.

If $D\ge4$, a subspace $E$ of $V_D$ is in $\ph'(V_D)$ if either
$E=0$ or if there exists
$E'\in V_{D-2}$ and $j\in[1,D+1]$ such that $E=\t_j(E')\op Fe_j$.
We then have
$$\th_j(\bx{E'})=\bx{E}\tag a$$
The following is proved in \cite{L23}.

(b) {\it The map $B\m<B>$ is a bijection $\ph(V_D)@>\si>>\ph'(V_D)$.}

For $B'\in\ph(V_{D-2})$ and $j\in[1,D+1]$ we define $t_j(B')\in\ph(V_D)$
by

(c) $<t_j(B')>=\t_j(<B'>)\op Fe_j$.
\nl
Using (a) we then have

(d) $\th_j<bx{<B'>}=\bx{<t_j(B')>}$.

For $x'\in V_{D-2},x\in V_D$ such that $B(x)=t_j(B'(x'))$ we show:

(e) $\th_j(\bx{x'})=\bx{x}+\bx{x+e_j}$.
\nl
Indeed, let $\e':\ph(V_{D-2})@>>>V_{D-2}$ be the bijection analogue to
$\e:\ph(V_D)@>>>V_D$. From \cite{L20, 1.9(b)} we have
$\e(B(x))=\t_j(\e'(B'(x'))+ce_j$ for some $c\in F$. In other words we
have $x=\t_j(x')+ce_j$. It follows that
$\bx{\t_j(x')}+\bx{\t_j(x')+e_j}=\bx{x}+\bx{x+e_j}$
so that (e) holds.

\subhead 2.3\endsubhead
For any $A\in\ph(V_D)$ we have

(a) $\bx{<A>}=\sum_{A'\in\ph(V_D)}d_A^{A'}\bx{\e\i(A')}$
\nl
where $d_A^{A'}=1$ if $\e\i(A')\in<A>$ and $d_A^{A'}=0$ if
$\e\i(A')\n<A>$;
in particular we have $d_A^{A'}\ne0\implies A'\le A$. Moreover we have
$d_A^A=1$.

It follows that for any $C\in\ph(V_D)$ we have

(b) $\bx{\e\i(C)}=\sum_{B\in\ph(V_D)}r_C^B\bx{<B>}$
\nl
with $r_C^B\in\ZZ$ uniquely determined and such that
$r_C^B\ne0\implies B\le C$; moreover we have $r_C^C=1$.

We deduce that for $A'\le A$ in $\ph(V_D)$ we have
$$d_A^{A'}=
\sum(-1)^kr_{A_0}^{A_1}r_{A_1}^{A_2}\do r_{A_{k-1}}^{A_k}\tag c$$
where the sum is taken over all sequences $A=A_0>A_1>A_2>\do>A_k=A'$ in
$\ph(V_D)$.

\proclaim{Theorem 2.4} If $C\in\ph(V_D),B\in\ph(V_D)$ satisfy $r_C^B\ne0$
then $|B|\le|C|$.
\endproclaim
We argue by induction on $D$. If $D=2$ the result is obvious.
We now assume that $D\ge4$. Let $x=\e\i(C)$.
We shall use a second induction on $\nu(x)$ (see 1.1).
If $\nu(x)=0$ then $x=0$ and $B=C=\emp$ so that $|C|=|B|=0$ and
the desired result holds. Now assume that $\nu(x)\ne0$ that is
$x\ne0$. Then for some $j\in[1,D+1]$ we have
$<C>=\t_j<C'>\op Fe_j$ for some $C'\in\ph(V_{D-2})$      
(In particular we have $\{j\}\in C$.)
Let $\e':\ph(V_{D-2})@>>>V_{D-2}$ be as in 2.2 and let $x'=\e'(C')$.
By the first induction hypothesis we have
$$\bx{x'}=\sum_{B'\in \ph(V_{D-2});|B'|\le|C'|}r'_{C'}{}^{B'}\bx{<B'>}$$
with $r'_{C'}{}^{B'}\in\ZZ$. Here the boxed entries refer to
$V_{D-2}$. Applying $\th_j$ (see 2.1) we deduce
$$\bx{\t_j(x')}+\bx{\t_j(x')+e_j}=
\sum_{B'\in\ph(V_{D-2});|B'|\le|C'|}r'_{C'}{}^{B'}\bx{<t_j(B')>}.$$
Using 2.2(e) this becomes
$$\bx{x}+\bx{x+e_j}=
\sum_{B'\in\ph(V_{D-2});|B'|\le|C'|}r'_{C'}{}^{B'}\bx{<t_j(B')>}.$$
Note that $|t_j(B')|=|B'|+1$. Similarly, $|C|=|C'|+1$. We see that

(a) $\bx{x}+\bx{x+e_j}$ is a linear combination of
elements $\bx{<B>}$ with $B\in\ph(V_D)$ such that $|B|\le|C|$.
\nl
Since $\{j\}\in C$, we see that $x+e_j\in<C>$ so that
$x+e_j<x$. This implies that $\nu(x+e_j)<\nu(x)$ (see 1.1).
By the second induction hypothesis we see that

(b) $\bx{x+e_j}$ is a linear combination of
elements $\bx{<B>}$ with $B\in\ph(V_D)$ such that $|B|\le|\e\i(x+e_j)|$.
\nl
Combining (a),(b), it remains to show that

(c) $|\e\i(x+e_j)|\le|C|$.
\nl
Since $\{\j\}\in C$, the object $C[j]\in\ph(V_D)$ is defined
as in \S1; it satisfies $\e(C[j])=\e(C)+e_j=x+e_j$, see 1.4. Thus,
$\e(C[j])=\e(\e\i(x+e_j))$. Using the fact that $\e:\ph(V_D)@>>>V_D$
is a bijection, we deduce that $C[j]=\e\i(x+e_j)$ so that
$|C[j]|=|\e\i(x+e_j)|$. From definitions we have
$|C[j]|=|C|-1$ (if $g_j(C)=1$) and
$|C[j]|=|C|$ (if $g_j(C)\ge2$).
In particular we have $|C[j]|\le|C|$ and 
(c) follows. The theorem is proved.

\proclaim{Corollary 2.5}If $A',A$ in $\ph(V_D)$ satisfy $A'\le A$,
then $|A'|\le|A|$.
\endproclaim
We can assume that $\e\i(A')\in<A>$ so that $d_A^{A'}=1$
(notation of 2.3). Using 2.3(c) we deduce that there exists
a sequence $A=A_0>A_1>A_2>\do>A_k=A'$ in $\ph(V_D)$ such that
$r_{A_0}^{A_1}\ne0,r_{A_1}^{A_2}\ne0,\do,r_{A_{k-1}}^{A_k}\ne0$.
Using 2.4, we see that
$$|A_0|\ge|A_1|\ge|A_2|\ge\do\ge|A_k|$$
so that $|A|\ge|A'|$, as desired.

\head 3. Matrix entries of the Fourier transform\endhead
\subhead 3.1\endsubhead
In this subsection $F$ (in 0.1)
is assumed to be a finite field with $q$
elements. Let $\ps:F@>>>\CC^*$ be a nontrivial homomorphism. Recall that
the Fourier transform $\cf_\ps:\CC[V_D]@>>>\CC[V_D]$ is the linear map
given by $\cf_\ps(\bx{x})=q^{-D/2}\sum_{y\in V_D}\ps((x,y))\bx{y}$
for all $x\in V_D$.

\subhead 3.2\endsubhead
We now return to the setup of 0.2.
Then $\ps$ in 3.1 is uniquely determined and
we write $\cf$ instead of $\cf_\ps$.
We are interested in the entries of the matrix of $\cf$ with respect to
the new basis, that is in the complex numbers $n_{B,B'}$ given by
$$\cf_D(\bx{<B>})=\sum_{B'\in\ph(V_D)}n_{B,B'}\bx{<B'>}$$
for any $B\in\ph(V_D)\}$.

We have the following result.

$$\text{If } n_{B,B'}\ne0\text{ and } B\ne B' \text{ then }
B<B' \text{ and }|B|<|B'|.\tag a$$
The fact that $|B|<|B'|$ is proved in \cite{L20a, 1.6}. A similar
argument shows that $B<B'$. Indeed, we can assume that $D\ge4$. As in
{\it loc.cit.}, it is enough to prove the following statement: if
$B_1<B'_1$ in $\ph(V_{D-2})$ and $j\in[1,D+1]$ then
$t_j(B_1)<t_j(B'_1)$ (compare \cite{L20, 1.20(b)}).

From (a) we see that the matrix $(n_{B,B'})$ is upper triangular. Its
diagonal entries are $\pm1$ since $\cf^2=1$.

For $B\in\ph(V_D)\}$ we have
$$\cf_D(\bx{<B>})=2^{|B|-(D/2)}\bx{<B>^\pe},$$
where $<B>^\pe=\{x\in V_D;(x,<B>)=0\}$ (since $<B>$ is an isotropic
subspace).
In particular we have
$$2^{-D/2}\bx{V_D}=\sum_{B'\in\ph(V_D)}n_{\emp,B'}\bx{<B'>}.$$

\subhead 3.3\endsubhead
For a sequence $i_1,i_2,\do,i_k$ in $[1,D+1]$ (identified with
$\ZZ/(D+1)$) we write 
$$[i_1i_2\do i_k]=\sum_J\bx{<\e\i(e_J)>}\in\CC[V_D]$$
where the sum is taken over all distinct subsets $J$
of $\ZZ/(D+1)$ of the form
$\{i_1+h,i_2+h,\do,i_k+h\}$ for some $h\in \ZZ/(D+1)$.

When $D\le8$, the element $\bx{V_D}\in\CC[V_D]$ can be written as follows.

$[1]-2[-]$  if $D=2$,

$[123]-4[-]$ if $D=4$,

$[1245]+[1235]+[1236]-2[123]-2[135]+4[13]-4[1]+8[-]$ if $D=6$,

$$\align&
[1246]+[123467]+[124567]-2[1234567]+2[123457]+2[134567]+2[123567]\\&
-2[13457]-2[12357]-2[13567]-4[12345]+4[1235]+4[1238]+4[147]\\&
-8[123]+8[13]-8[14]+16[-]\endalign$$
if $D=8$.

Here $-$ denotes the empty sequence.

We state:
\proclaim{Conjecture 3.4} For any $B,B'$ we have either
$n_{B,B'}=0$ or $n_{B,B'}=\pm2^{-t}$ for some $t\in[0,D/2]$.
\endproclaim
One can show that for $B,B'$ in $\ph(V_D)$, the coefficient $n_{B,B'}$
is either $0$ or is equal to some $n_{\emp,B'_1}$ where
$B'_1\in\ph(<B>^\pe/<B>)$ (note that by \cite{L22}, $<B>^\pe/<B>$
inherits
from $V_D$ a structure similar to that of $V_D$).
In this way we see that it is enough to prove the conjecture in the
special case where $B=\emp$. (This argument is similar to that in
\cite{L20a, 1.5,1.6}.) Using this and the results in 3.2 we see that the
conjecture holds for $D\le8$.

\subhead 3.5\endsubhead
We describe the matrix $(n_{B,B'})$ assuming that $D=2$ or $D=4$.

If $D=2$ the matrix is
$$\left(\mat -1 &1/2&1/2&1/2\\
              0 &1 &0&0     \\
              0 &0 &1&0     \\
              0 &0 &0&1 \endmat\right)$$
If $D=4$ the matrix is
                               
$$\left(\mat -1 &  0&0&0&0&0  &0&0&0&0&0  & 1/4&1/4&1/4&1/4&1/4\\
              0 & -1&0&0&0&0& 0&0&1/2&1/2&0&   1/2&0&0&0&0\\
              0 & 0&-1&0&0&0& 1/2&0&0&0&1/2&   0&1/2&0&0&0\\
              0 & 0&0&-1&0&0& 0&1/2&1/2&0&0&   0&0&1/2&0&0\\
              0 & 0&0&0&-1&0& 0&0&0&1/2&1/2&   0&0&0&1/2&0\\
              0 & 0&0&0&0&-1& 1/2&1/2&0&0&0&   0&0&0&0&1/2\\
              0 & 0&0&0&0&0& 1&0&0&0&0&   0&0&0&0&0\\
              0 & 0&0&0&0&0& 0&1&0&0&0&   0&0&0&0&0\\
              0 & 0&0&0&0&0& 0&0&1&0&0&   0&0&0&0&0\\
              0 & 0&0&0&0&0& 0&0&0&1&0&   0&0&0&0&0\\
              0 & 0&0&0&0&0& 0&0&0&0&1&   0&0&0&0&0\\
              0 & 0&0&0&0&0& 0&0&0&0&0&   1&0&0&0&0\\
              0 & 0&0&0&0&0& 0&0&0&0&0&   0&1&0&0&0\\
             0 & 0&0&0&0&0& 0&0&0&0&0&   0&0&1&0&0\\
             0 & 0&0&0&0&0& 0&0&0&0&0&   0&0&0&1&0\\
             0 & 0&0&0&0&0& 0&0&0&0&0&   0&0&0&0&1 \endmat\right)$$

\subhead 3.6\endsubhead
The examples in 3.3 suggest that when $n_{\emp,B'}\ne0$ then the number
$t$ in 3.4 is $(D/2)-k$ where $k$ is largest number such that
$B'=B'_0<B'_1<\do<B'_k$ for some $B'_1,\do,B'_k$ in $\ph(V_D)$.
It would be interesting to find a criterion for when $n_{\emp,B'}$
is $<0,0$ or $>0$.

\subhead 3.7\endsubhead
In the setup of 3.1, it would be interesting to extend the definition
of the new basis of $\CC[V_D]$ to the case where $q>2$ in such a way
that the Fourier
transform in 3.1 is again triangular with respect to this basis.

\enddocument